\documentclass[11pt,reqno,a4paper]{amsart}
\usepackage[a4paper,margin=30mm,headheight=14pt,headsep=7mm]{geometry}
\usepackage[T1]{fontenc}
\usepackage{lmodern}
\usepackage{amsmath,amssymb,mathtools}
\usepackage[expansion=false]{microtype}
\usepackage{needspace}
\usepackage[hidelinks]{hyperref}
\usepackage[nobysame,initials]{amsrefs}
\hypersetup{pdftitle={Hereditary lowerability of diagonal entropy scales},
 pdfauthor={Xiaojun Huang and Zhuang Li},
 pdfsubject={Upper capacity entropy and compact subset realization},
 pdfkeywords={entropy scale, hereditary lowerability, entropy dimension},
 pdflang={en-US},pdfdisplaydoctitle=true,bookmarksnumbered=true}
\numberwithin{equation}{section}
\newtheorem{theorem}{Theorem}[section]
\newtheorem{proposition}[theorem]{Proposition}
\newtheorem{lemma}[theorem]{Lemma}

\theoremstyle{definition}
\newtheorem{definition}[theorem]{Definition}

\theoremstyle{remark}
\newtheorem{remark}[theorem]{Remark}
\newcommand{\N}{\mathbb N}
\newcommand{\Nzero}{\mathbb N_0}
\newcommand{\Z}{\mathbb Z}
\newcommand{\scl}{\operatorname{scl}}

\newcommand{\card}[1]{\lvert #1\rvert}
\newcommand{\Dp}{D_p}

\newcommand{\Exp}{\mathcal E}

\newcommand{\eps}{\varepsilon}

\title[Hereditary lowerability of diagonal entropy scales]{Hereditary lowerability of diagonal entropy scales}
\author{Xiaojun Huang}
\address{College of Mathematics and Statistics, Chongqing University, Chongqing 401331, China}
\email{hxj@cqu.edu.cn}
\author{Zhuang Li}
\address{College of Mathematics and Statistics, Chongqing University, Chongqing 401331, China}
\email{jiu\_L@qq.com}
\date{}
\begin{document}
\raggedbottom
\setlength{\emergencystretch}{1em}
\begin{abstract}
Let $T$ be a continuous self-map of a compact metric space. For each integer
$p\geq2$, we consider the upper capacity entropy scale obtained by dividing
the $p$-fold logarithm of an orbit covering number by the $(p-1)$-fold
logarithm of time. We prove that every nonempty compact set contains a
compact subset realizing any value between zero and its own scale entropy.
Every strictly intermediate value is realized by a countable compact set
with exactly one accumulation point. The proof combines localization at a
fixed spatial resolution, a finite-tree selection lemma, and refinements
of orbit covers at dyadic times. An example shows that the assertion about
one accumulation point does not extend to the upper endpoint. On the product-diagonal construction used in classical entropy
lowerability, the diagonal scales preserve the values of the base subsets,
although classical upper capacity entropy takes only zero or infinity.
\end{abstract}
\maketitle

\section{Introduction}

The entropy of a compact subset describes the growth of distinguishable
orbit segments starting in that set. The set need not be invariant. This
leads to an intermediate-value question: if a compact set has a given
entropy, which smaller values occur on its compact subsets? The answer
depends both on the growth normalization and on the definition of subset
entropy.

The open-cover definition of topological entropy was introduced by
Adler, Konheim, and McAndrew \cite{AKM}. Bowen's metric formulation
\cite{Bowen71} describes the same invariant on compact systems through
orbit separation. For arbitrary subsets one can retain a common orbit
length in every covering ball, obtaining upper capacity entropy, or allow
the orbit lengths to vary, as in Bowen's definition for noncompact sets
\cite{Bowen73}. These constructions can assign different values even to
compact sets that are not invariant.

Huang, Ye, and Zhang \cite{HYZ10} formulated lowerability and hereditary
lowerability for upper capacity entropy. Lowerability asks for all values
between zero and the entropy of the whole system; hereditary lowerability
asks the same question inside each nonempty compact subset. They proved
lowerability for systems of finite topological entropy. Their subsequent
work \cite{HYZ14} established lowerability without the finiteness assumption
and constructed systems that are not hereditarily lowerable. In particular,
their product-diagonal example contains a compact set whose compact
subsets have classical upper capacity entropy only zero or infinity. In
a recent preprint, Wang \cite{Wang} proves hereditary lowerability for
finite-entropy homeomorphisms. Thus the finite-entropy and infinite-entropy
settings have different behavior for the classical normalization.

To distinguish orbit growth below the exponential scale, one can
replace the classical normalization by slower growth functions. Katok and Thouvenot \cite{KT} introduced slow-entropy invariants in
measure-preserving dynamics. In topological dynamics, de Carvalho
\cite{Carvalho} and Dou, Huang, and Park \cite{DHP} developed entropy
dimension using growth rates between
polynomial and exponential ones. More generally, Helfter \cite{Helfter}
organized metric covering invariants by families of gauge functions,
including families obtained from iterated exponentials and logarithms.
Xie, Chen, and Zhou \cite{XCZ} formulate
Bowen and upper capacity entropy scales in dynamics and study their
comparison and variational properties. These results motivate asking how
lowerability behaves for a specified family of orbit-growth scales.

Here the time convention is fixed as $t=e^{-n}$. If $r_n(E,\eps)$ denotes
the least number of radius-$\eps$ Bowen balls of length $n$ covering $E$,
the diagonal member of the resulting hierarchy is
\begin{equation}\label{eq:intro-D}
 \Dp(T,E)=\lim_{\eps\downarrow0}\limsup_{n\to\infty}
 \frac{\ell_p(r_n(E,\eps))}{L_{p-1}(n)},\qquad p\geq2.
\end{equation}
Here $L_k$ is the $k$-fold logarithm, defined for sufficiently large
arguments, and $\ell_k$ is a regularized version defined on $[0,\infty)$;
both are specified in Section~\ref{sec:prelim}. The regularization makes
\eqref{eq:intro-D} meaningful when covering numbers stay bounded. For
$p=2$ and $E=X$, this is upper entropy dimension; for general $E$ it
is the corresponding upper capacity extension to subsets. At every
fixed resolution, orbit covering numbers are
at most exponential in $n$, so $0\leq\Dp(T,E)\leq1$.

\begin{theorem}\label{thm:main}
Let $X$ be a nonempty compact metric space and let $T\colon X\to X$ be
continuous. Fix an integer $p\geq2$. For every nonempty compact set
$E\subset X$ and every $a\in[0,\Dp(T,E)]$, there is a nonempty compact set
$K_a\subset E$ such that
\[
 \Dp(T,K_a)=a.
\]
If $0<a<\Dp(T,E)$, the set $K_a$ can be chosen countable with exactly one
accumulation point.
\end{theorem}

Consequently every compact continuous system is hereditarily lowerable
for each fixed diagonal scale in \eqref{eq:intro-D}. Neither $E$ nor $K_a$
is required to be invariant. No finiteness assumption on classical
topological entropy is imposed, and invertibility is not used. The
statement is separate for each $p$; it does not prescribe several scale
values on the same subset.

The construction follows the extraction strategy for countable compact sets
in \cite{Wang}*{Lemma~3.5}, using metric covering trees and the diagonal
scale estimates developed below. First, fixed-resolution localization
produces a countable compact source set with one accumulation point and
scale entropy larger than the target. Second, finite subsets of successively
smaller tails are selected so that their orbit counts satisfy one common
upper bound at all times. For this step we use Wang's finite-tree lemma
\cite{Wang}*{Lemma~3.3}. The trees here come from successive metric covers
at dyadic times. Their node counts include a logarithmic loss, which the
parameter gap in the diagonal scales absorbs. A recurring schedule of
spatial resolutions supplies the matching lower bound.

The restriction to strictly intermediate values in the assertion about
accumulation points is necessary. Section~\ref{sec:supplements} gives,
for every fixed $p\geq2$ and $H\in(0,1]$, a compact set with $\Dp=H$
for which every compact subset with at most one accumulation point
has $\Dp<H$. The same section computes $\Dp$ on the
product-diagonal example of \cite{HYZ14}: the classical zero-or-infinity
restriction coexists with preservation of the base subsets' diagonal
scale values.

Section~\ref{sec:prelim} records the definitions, covering facts, and
finite combinatorial terminology used later. Section~\ref{sec:calculus}
establishes the scale estimates. Sections~\ref{sec:local} and
\ref{sec:extraction} give the localization and the proof of
Theorem~\ref{thm:main}.

\section{Preliminaries}\label{sec:prelim}

\subsection{Orbit metrics and subset entropy}
Let $(X,d)$ be a nonempty compact metric space and
let $T\colon X\to X$ be continuous. We use $\N=\{1,2,\ldots\}$ and
$\Nzero=\{0,1,\ldots\}$, and all logarithms are natural unless a base is
displayed. For $n\in\N$, set
\[
 d_n(x,y)=\max_{0\leq j<n}d(T^jx,T^jy),\qquad
 B_n(x,\eps)=\{y\in X:d_n(x,y)<\eps\}.
\]
For a nonempty set $E\subset X$, let $r_n(E,\eps)$ be the least number of
such balls covering $E$, with centers allowed in $X$. Let $s_n(E,\eps)$
be the largest cardinality of a subset $F\subset E$ satisfying
$d_n(x,y)>\eps$ for distinct $x,y\in F$. Such a set is called
$(n,\eps)$-separated. The metric $d_n$ induces the topology of $X$, so
compactness makes these counts finite. The largest separated cardinality
exists even if $E$ is not closed: possible cardinalities are integers in
a nonempty finite set. Both counts are nondecreasing in $n$.

The standard covering--separation inequalities, in the present open-ball
convention, are
\begin{equation}\label{eq:rs}
 r_n(E,2\eps)\leq s_n(E,\eps)\leq r_n(E,\eps/2).
\end{equation}
Indeed, a separated set of largest cardinality is maximal under inclusion,
so every point of $E$ is within $d_n$-distance at most $\eps$ of it. Also,
an open ball of radius $\eps/2$ contains at most one point of an
$(n,\eps)$-separated set. See \cite{HYZ14}*{Section~2} for the customary
covering and separation descriptions of subset entropy.

Classical upper capacity entropy is
\begin{equation}\label{eq:classical-entropy}
 h(T,E)=\lim_{\eps\downarrow0}\limsup_{n\to\infty}
 \frac{\log r_n(E,\eps)}{n}.
\end{equation}
For $E=X$ this is topological entropy. The counts in
\eqref{eq:classical-entropy} use one common orbit length at each stage.
Bowen's entropy for noncompact sets instead uses weighted covers whose
orbit lengths may vary from ball to ball \cite{Bowen73}. All lowerability
assertions below concern upper capacity counts.

\subsection{Iterated functions and the time convention}
Set
\[
 L_0(u)=u,\quad L_{k+1}(u)=\log L_k(u),\qquad
 \Exp_0(u)=u,\quad \Exp_{k+1}(u)=\exp(\Exp_k(u)).
\]
An expression containing $L_k$ is used only where its successive
logarithms are defined, and in limits its argument is sufficiently large
that they are positive. Thus $L_1(u)=\log u$, $L_2(u)=\log\log u$, and
$\Exp_2(u)=\exp(e^u)$. The identities
$L_k(\Exp_k(u))=u$ and $\Exp_k(L_k(u))=u$ hold on these domains.
For counts that may equal one, use
\begin{equation}\label{eq:regularized}
 \ell_0(u)=u,\qquad \ell_{k+1}(u)=\log(e+\ell_k(u)),\quad u\geq0.
\end{equation}
Every $\ell_k$ is increasing, and
\begin{equation}\label{eq:reg-asymp}
 \ell_k(u)=L_k(u)+o(1)\qquad(u\to\infty),\quad k\geq1.
\end{equation}
For $k=1$, the difference is $\log(1+e/u)$. If
$\ell_k(u)=L_k(u)+o(1)$, the next difference is
$\log(1+(e+o(1))/L_k(u))$, proving \eqref{eq:reg-asymp} by induction.

In asymptotic notation, $F=O(G)$ means $|F|\leq C|G|$ eventually for some
constant $C$, and $F=o(G)$ means $F/G\to0$. A subscript on $O$, such as
$O_{\eps}(1)$, permits the constant to depend on that subscript. Limits
in time are always taken before the spatial resolution tends to zero.

For positive integers $p,q$ and $\alpha>0$, the iterated scale family of
\cite{Helfter}*{Proposition~2.4} can, for sufficiently small $t>0$, be
written as
\begin{equation}\label{eq:scaling}
 \scl^{p,q}_\alpha(t)=
 \frac{1}{\Exp_p(\alpha L_q(t^{-1}))}.
\end{equation}
Only its behavior as $t\downarrow0$ is needed. Under the convention
$t=e^{-n}$, its reciprocal becomes
\begin{equation}\label{eq:growth-family}
 A^{p,q}_\alpha(n)=\Exp_p(\alpha L_{q-1}(n)).
\end{equation}
We define
\begin{equation}\label{eq:critical-def}
 h^{p,q}(T,E)=\lim_{\eps\downarrow0}\inf\left\{
 \alpha>0:\frac{r_n(E,\eps)}{A^{p,q}_\alpha(n)}\longrightarrow0
 \right\}.
\end{equation}
Here $\inf\varnothing=+\infty$. This explicitly fixes the normalization;
changing the substitution for
$t$ would give a different time-growth question. The notation used in the
rest of the paper is $\Dp(T,E)=h^{p,p}(T,E)$ for $p\geq2$.

\begin{lemma}\label{lem:critical}
For every nonempty $E\subset X$ and all $p,q\geq1$,
\begin{equation}\label{eq:log-formula}
 h^{p,q}(T,E)=\lim_{\eps\downarrow0}\limsup_{n\to\infty}
 \frac{\ell_p(r_n(E,\eps))}{L_{q-1}(n)}.
\end{equation}
The same outer limit is obtained with $s_n(E,\eps)$ in place of
$r_n(E,\eps)$. These quantities are independent of the compatible metric
on $X$. In particular, $h^{1,1}(T,E)=h(T,E)$.
\end{lemma}

\begin{proof}
Fix $\eps$ and put
\[
 R_n=r_n(E,\eps),\qquad v_n=L_{q-1}(n),\qquad
 \beta=\limsup_n\frac{\ell_p(R_n)}{v_n}.
\]
If $\alpha>\beta$, choose
$\beta<b'<b<\alpha$. Eventually
$\ell_p(R_n)<b'v_n<\ell_p(\Exp_p(bv_n))$, by
\eqref{eq:reg-asymp}. Monotonicity gives
$R_n\leq\Exp_p(bv_n)$, and
\[
 \frac{\Exp_p(bv_n)}{\Exp_p(\alpha v_n)}\longrightarrow0.
\]
For $p=1$ this is $e^{-(\alpha-b)v_n}$; for $p>1$ its logarithm is
$\Exp_{p-1}(bv_n)-\Exp_{p-1}(\alpha v_n)\to-\infty$.
Conversely, if $0<\alpha<\beta$, choose
$\alpha<b<b'<\beta$. On an infinite subsequence,
$\ell_p(R_n)>b'v_n>\ell_p(\Exp_p(bv_n))$. Hence
$R_n\geq\Exp_p(bv_n)$ on that subsequence, and the ratio with
$\Exp_p(\alpha v_n)$ tends to infinity there. These two arguments also
cover $\beta=0$ and $\beta=+\infty$ by the available one-sided
inequality. Bounded sequences $R_n$ are included in the first argument.
This proves the critical-parameter formula at fixed $\eps$.

Inequalities \eqref{eq:rs} give the separated-count version after
$\eps\downarrow0$. If $d'$ is another compatible metric, compactness
gives, for every $\eps>0$, a $\delta>0$ such that
$d(x,y)<\delta$ implies $d'(x,y)<\eps$. The same implication holds for
the corresponding orbit metrics, whence
$r_n(E,d',\eps)\leq r_n(E,d,\delta)$. Interchanging the two metrics
proves metric independence. Finally,
$0\leq\ell_1(R)-\log R\leq\log(1+e)$ for $R\geq1$; division by $n$
gives the last assertion.
\end{proof}

For $b>0$, two examples of numerical growth rates are
\[
 R_n=\exp(n^b),\qquad \frac{\ell_2(R_n)}{\log n}\longrightarrow b,
\]
and
\[
 R_n=\exp(\exp((\log n)^b)),\qquad
 \frac{\ell_3(R_n)}{\log\log n}\longrightarrow b.
\]
Bounded counts have rate zero. Lemma~\ref{lem:range} bounds the rates
that can occur dynamically.

In the case $p=q=2$, one also has
\begin{equation}\label{eq:entropy-dimension}
 D_2(T,E)=\lim_{\eps\downarrow0}\inf\left\{\alpha>0:
 \limsup_{n\to\infty}\frac{\log r_n(E,\eps)}{n^\alpha}=0\right\}.
\end{equation}
To check this, fix $\eps$ and let $\beta$ be the limsup in
\eqref{eq:log-formula} for $p=q=2$. If $\beta<b<\alpha$, then
$\log r_n(E,\eps)\leq n^b$ eventually. If $\alpha<b<\beta$, then
$\log r_n(E,\eps)\geq n^b/2$ on an infinite subsequence.
These assertions follow by undoing the two regularized logarithms;
the bounded difference between $\ell_1(R)$ and $\log R$ does not
change either estimate. Division by $n^\alpha$ proves
\eqref{eq:entropy-dimension}, with the same one-sided argument at
$\beta=0$. For $E=X$, this is the upper entropy dimension considered
in \cites{Carvalho,DHP}; equation \eqref{eq:entropy-dimension} specifies
the extension used here when $E$ is not invariant.

For the estimates used later, fix $p\geq2$, put $r=p-1$, and define,
for $b>0$ and sufficiently large real $t$,
\begin{equation}\label{eq:fG}
 f_b(t)=\Exp_r(bL_r(t)),\qquad G_b(t)=\exp(f_b(t)).
\end{equation}
Thus $G_b=A^{p,p}_b$, and $f_1(t)=t$. For $p=2$ one has
$f_b(t)=t^b$ and $G_b(t)=e^{t^b}$; for $p=3$ one has
$f_b(t)=\exp((\log t)^b)$. Thresholds will always be chosen with
$L_r(t)>0$, on which domain both functions increase in $t$ and $b$.

\subsection{Finite unions, localization rates, and the range of values}
Fix $p\geq2$. To retain a spatial resolution in the notation, put
\begin{equation}\label{eq:fixed-rate}
 d_\eps(E)=\limsup_{n\to\infty}
 \frac{\ell_p(s_n(E,\eps))}{L_{p-1}(n)}.
\end{equation}
Thus $\Dp(T,E)=\lim_{\eps\downarrow0}d_\eps(E)$.
The symbol $d_\eps(E)$ denotes a rate, whereas $d_n(x,y)$ denotes an
orbit distance. The map $T$ and the fixed integer $p$ are suppressed in
\eqref{eq:fixed-rate}.

\begin{lemma}\label{lem:finite}
The quantities $d_\eps(E)$ and $\Dp(T,E)$ are monotone under inclusion.
For a finite family of nonempty subsets $E_1,\ldots,E_m$ of $X$,
\begin{align}
 d_\eps\left(\bigcup_{i=1}^m E_i\right)&=\max_i d_\eps(E_i),\label{eq:finite-fixed}\\
 \Dp\left(T,\bigcup_{i=1}^m E_i\right)&=\max_i\Dp(T,E_i).\label{eq:finite-full}
\end{align}
Finite nonempty sets have value zero. Deleting finitely many points from
an infinite set changes neither rate.
\end{lemma}

\begin{proof}
Let $S_n=\max_i s_n(E_i,\eps)$. Assign each point of a separated set
in the union to one containing $E_i$. The assigned points in each $E_i$
are separated, so
\[
 S_n\leq s_n\left(\bigcup_iE_i,\eps\right)
 \leq\sum_i s_n(E_i,\eps)\leq mS_n.
\]
For fixed $m$ and $p$, the difference
$\ell_p(mu)-\ell_p(u)$ is bounded on $u\geq1$. For $p=1$ it is at most
$\log m$, and induction follows because $u\mapsto\log(e+u)$ has
derivative at most $1/e$. Division by $L_{p-1}(n)\to\infty$ removes
this bounded difference. The limsup of a maximum of finitely many
sequences is the maximum of their limsups, proving
\eqref{eq:finite-fixed}. Taking the monotone resolution limit proves
\eqref{eq:finite-full}. A finite set has uniformly bounded counts.
Writing a set as the union of its remaining points and its deleted
finite part proves the deletion statement.
\end{proof}

\begin{lemma}\label{lem:range}
For every $\eps>0$ there is an integer $M_\eps\geq1$ such that
\begin{equation}\label{eq:exponential-cover}
 r_n(E,\eps)\leq r_n(X,\eps)\leq M_\eps^n
 \qquad(n\in\N,\ \varnothing\ne E\subset X).
\end{equation}
Consequently $0\leq\Dp(T,E)\leq1$ for $p\geq2$, and
$h^{p,q}(T,E)=0$ if $p>q$. If $h(T,E)>0$, then $\Dp(T,E)=1$ for every
$p\geq2$.
\end{lemma}

\begin{proof}
Cover $X$ by $M_\eps$ open sets of diameter less than $\eps$. Each
nonempty intersection
$U_{i_0}\cap T^{-1}U_{i_1}\cap\cdots\cap T^{-(n-1)}U_{i_{n-1}}$
has $d_n$-diameter less than $\eps$. There are at most $M_\eps^n$
such intersections, and one point from each supplies an open Bowen-ball
cover. This proves \eqref{eq:exponential-cover}.
For $M>1$ and $p\geq2$,
\[
 \frac{\ell_p(M^n)}{L_{p-1}(n)}\longrightarrow1.
\]
If $M=1$ the limit is zero. This gives the bound on $\Dp$; when $p>q$,
use additionally $L_{p-1}(n)/L_{q-1}(n)\to0$.
If $h(T,E)>0$, some fixed $\eps$ and $c>0$ satisfy
$r_n(E,\eps)\geq e^{cn}$ for infinitely many $n$. The corresponding
normalized $p$-fold logarithms tend to one on that subsequence.
The upper bound just established gives equality.
\end{proof}

\begin{definition}\label{def:lowerable}
For a fixed $p\geq2$, $(X,T)$ is \emph{lowerable for $\Dp$} if each
$a\in[0,\Dp(T,X)]$ occurs on a nonempty compact subset of $X$. It is
\emph{hereditarily lowerable for $\Dp$} if the same assertion holds
inside every nonempty compact subset $E$, with upper endpoint
$\Dp(T,E)$.
\end{definition}
These definitions follow the subset formulation of \cites{HYZ10,HYZ14}.
An accumulation point of a set is a point every neighborhood of which
contains a different point of that set. A model for the compact sets used
below is $\{0\}\cup\{1/j:j\in\N\}\subset\mathbb R$: it is countable,
and zero is its only accumulation point. On a compact metric space an
infinite set has an accumulation point, so a compact set with none is
finite.

\subsection{Finite partitions and rooted trees}
A finite partition of a finite set $P$ is a collection of nonempty,
pairwise disjoint subsets, called atoms, whose union is $P$. The join
$\mathcal P\vee\mathcal Q$ consists of all nonempty intersections of
an atom of $\mathcal P$ with an atom of $\mathcal Q$. It refines both
partitions: each new atom is contained in one atom of either partition.
For instance, the join of
$\{\{1,2\},\{3,4\}\}$ and $\{\{1,3\},\{2,4\}\}$ form the singleton
partition of $\{1,2,3,4\}$.

A finite rooted tree is a finite connected graph with no cycles and
with a distinguished vertex, its root. A vertex has level $k$ if the
path from it to the root has $k$ edges. Its parent is the preceding
vertex on that path. A leaf is a vertex with no children. We consider
trees whose leaves all lie at the same level $J$, and call these
vertices terminal leaves. A vertex is occupied by a chosen set of leaves if at least one
chosen leaf descends from it. For example, a binary tree of depth two
has leaves $00,01,10,11$. Choosing $00$ and $01$ occupies one vertex at
level one and two at level two; choosing $00$ and $10$ occupies two
vertices at each of these levels.

A sequence of refining partitions of $P$ gives such a tree: atoms are
vertices, inclusion determines parenthood, and singleton atoms in the
last partition are leaves. Controlling the number of occupied vertices
therefore controls how many partition atoms are needed to cover a
selected subset of $P$.

\begin{lemma}[Wang, {\cite{Wang}*{Lemma~3.3}}]\label{lem:tree-selection}
Let a finite rooted tree have $A$ leaves at a common level $J$, and let
$B_k$ be the number of vertices at level $k$ having a terminal descendant.
For every integer $1\leq M\leq A$, one can select exactly $M$ leaves
such that the number $b_k$ of occupied vertices satisfies
\begin{equation}\label{eq:tree-bound}
 b_k\leq2k^2\left(1+\frac{MB_k}{A}\right),\qquad 1\leq k\leq J.
\end{equation}
\end{lemma}
We use this finite combinatorial result in the stated form. Its proof in
\cite{Wang} applies to arbitrary finite rooted trees and requires no
dynamical assumptions. In particular, its loss depends on the current
level $k$, not on the final depth $J$.

\subsection{Symbolic and product models}
The two-sided full shift is
$\Sigma_2=\{0,1\}^{\Z}$ with $(\sigma x)_i=x_{i+1}$. A compatible
metric bounded by one is
\begin{equation}\label{eq:shift-metric}
 d_0(x,y)=\frac13\sum_{i\in\Z}2^{-|i|}|x_i-y_i|.
\end{equation}
A cylinder set specifies the coordinates in a finite index set and
leaves the others unrestricted. Cylinders are open and closed and form
a basis for the product topology. Distinct patterns on $[0,n)$ have
Bowen distance at least $1/3$, whereas agreement on $[-q,n-1+q]$
gives Bowen distance at most $(2/3)2^{-q}$. Thus, for
$0<\eps<1/6$ and $(2/3)2^{-q}<\eps$,
\[
 2^n\leq r_n(\Sigma_2,\eps)\leq2^{n+2q}.
\]
The lower bound uses one representative for each pattern on $[0,n)$;
the upper bound uses cylinders on $[-q,n-1+q]$. Hence
$h(\sigma,\Sigma_2)=\log2$, and Lemma~\ref{lem:range} gives
$\Dp(\sigma,\Sigma_2)=1$ for every $p\geq2$.

For any compact metric space $(X,d)$ with $d\leq1$, the countable product
$X^{\N}$ is compact and metrizable by
\begin{equation}\label{eq:product-metric}
 \rho(y,z)=\sum_{j\geq1}2^{-j}d(y_j,z_j).
\end{equation}
For completeness, every sequence in the product has a subsequence
converging in coordinate one, a further subsequence converging in
coordinate two, and so on. The diagonal subsequence converges in every
coordinate. To obtain convergence for \eqref{eq:product-metric}, first
make the tail sum smaller than a prescribed tolerance, and then use
convergence of the finitely many remaining coordinates. Sequential
compactness is compactness for a metric space. This also explains why
only finitely many coordinates are visible at a fixed resolution.

\section{Estimates for the diagonal scales}\label{sec:calculus}

Fix $p\geq2$, write $r=p-1$, and use the functions $f_b,G_b$
defined in \eqref{eq:fG}. The following estimates account for the
changes in covering counts introduced by the tree construction.

\Needspace{7\baselineskip}
\begin{lemma}\label{lem:scale-calculus}
For fixed positive parameters the following statements hold.
\begin{enumerate}
\item For $b>0$,
\begin{equation}\label{eq:scale-basic}
 \frac{\ell_p(G_b(t))}{L_r(t)}\longrightarrow b,
 \qquad \frac{f_b(t)}{\log t}\longrightarrow+\infty.
\end{equation}
\item If $0<b<d$, then
\begin{equation}\label{eq:absorb}
 (\log t)f_b(t)=o(f_d(t)).
\end{equation}
\item For $b>0$, $A\geq0$, and $C,c>0$,
\begin{equation}\label{eq:polynomial-loss}
 \frac{\ell_p(Ct^A G_b(ct))}{L_r(t)}\longrightarrow b.
\end{equation}
The same conclusion holds if a fixed nonnegative constant is added
inside $\ell_p$.
\item For $u,v>0$, $f_u\circ f_v=f_{uv}$ wherever both sides are
defined. In particular, the inverse of $f_D$ on a sufficiently large
interval is $f_{1/D}$.
\item Suppose $0<c<D$ and $0<a<D$. Set
\begin{equation}\label{eq:threshold-model}
 u_N=f_D^{-1}\left(\tfrac12 f_c(N)\right),\qquad
 M_N=\lfloor G_a(u_N)\rfloor.
\end{equation}
Then $u_N\to\infty$, $u_N<N$ eventually, and
\begin{equation}\label{eq:threshold-limit}
 \frac{\ell_p(M_N)}{L_r(N)}\longrightarrow\frac{ac}{D}.
\end{equation}
\end{enumerate}
\end{lemma}

\begin{proof}
We first record two comparisons of iterated exponentials. For $d>b>0$
and $j\geq1$,
\begin{equation}\label{eq:tower-comparisons}
 \frac{\Exp_j(dz)}{\Exp_j(bz)}\longrightarrow\infty,
 \qquad
 \frac{\Exp_j(bz)}{\Exp_{j-1}(z)}\longrightarrow\infty
 \quad(z\to\infty).
\end{equation}
For $j=1$, these are $e^{(d-b)z}\to\infty$ and
$e^{bz}/z\to\infty$. If either comparison holds at $j$, the logarithm
of the corresponding ratio at $j+1$ is the difference between its two
level-$j$ terms. The larger term divided by the smaller one tends to
infinity, and both terms tend to infinity; hence their difference tends
to infinity. This proves both assertions by induction.

The identity $L_p(G_b(t))=bL_r(t)$ and
\eqref{eq:reg-asymp} give the first limit in
\eqref{eq:scale-basic}. For $r=1$, the second limit and
\eqref{eq:absorb} reduce to
$t^b/\log t\to\infty$ and $(\log t)t^{b-d}\to0$.
If $r\geq2$, put $z=L_r(t)$. Taking logarithms gives
\begin{align*}
 \log\frac{f_b(t)}{\log t}
 &=\Exp_{r-1}(bz)-\Exp_{r-2}(z),\\
 \log\frac{f_d(t)}{(\log t)f_b(t)}
 &=\Exp_{r-1}(dz)-\Exp_{r-1}(bz)-\Exp_{r-2}(z).
\end{align*}
Both right-hand sides tend to infinity by
\eqref{eq:tower-comparisons}: in the second line each subtracted term
is of smaller order than the first one. This proves the first two parts.

For part (3), first note that $L_r(ct)/L_r(t)\to1$. This follows from
$\log(ct)=\log t+\log c$ for $r=1$, and from repeated logarithms for
larger $r$. Part (1) gives
$\log(Ct^A)=O(\log t)=o(f_b(ct))$, whence
\[
 \log(Ct^A G_b(ct))=f_b(ct)(1+o(1)).
\]
Taking $r$ further logarithms changes $L_r(f_b(ct))$ by $o(1)$.
That expression equals $bL_r(ct)$. Equation
\eqref{eq:reg-asymp} now proves \eqref{eq:polynomial-loss}.
Adding a fixed constant changes the first logarithm by $o(1)$ and
therefore gives the same limit.

The cancellation of $L_r$ with $\Exp_r$ gives
\[
 f_u(f_v(t))=\Exp_r\bigl(uL_r(\Exp_r(vL_r(t)))\bigr)
 =\Exp_r(uvL_r(t)),
\]
which proves part (4).
For part (5), $f_c(N)\to\infty$, so $u_N\to\infty$. Since $c<D$
and $L_r(N)>0$, one has $\tfrac12 f_c(N)<f_D(N)$ eventually;
strict monotonicity gives $u_N<N$. Moreover,
\begin{equation}\label{eq:half-threshold}
 D L_r(u_N)=L_r\left(\tfrac12 f_c(N)\right)
 =cL_r(N)+O(1).
\end{equation}
For $r=1$ the error is exactly $-\log2$, and for $r\geq2$ it tends
to zero. Since $G_a(u_N)\to\infty$, taking the integer part changes
its first logarithm by $o(1)$ and hence its $p$-fold logarithm by
$o(1)$. Thus
$\ell_p(M_N)=aL_r(u_N)+o(1)$, and
\eqref{eq:half-threshold} proves \eqref{eq:threshold-limit}.
\end{proof}

In the construction below, the parameters in
Lemma~\ref{lem:scale-calculus} are fixed separately at each stage before
the time is chosen large enough for the required error bound. No
uniformity over parameters from different stages is needed.

\section{Localization and metric covering trees}\label{sec:local}

Keep $p\geq2$ and $r=p-1$ fixed.

\subsection{A countable compact source}

The next two lemmas adapt the localization approach of
\cite{HYZ10}*{Section~5} to the fixed-resolution rates
\eqref{eq:fixed-rate} used in the extraction below.

\begin{lemma}\label{lem:local-point}
Let $E\subset X$ be nonempty and compact. If $d_\eps(E)>b>0$, then
there is $x_*\in E$ such that every relative open neighborhood $U$ of
$x_*$ in $E$ satisfies $d_\eps(U)\geq b$.
\end{lemma}

\begin{proof}
Otherwise, each $x\in E$ has a relative open neighborhood $U_x$ with
$d_\eps(U_x)<b$. Choose a finite subcover
$U_{x_1},\ldots,U_{x_m}$ of $E$. Since these sets are contained in
$E$, their union is $E$, and \eqref{eq:finite-fixed} gives
$d_\eps(E)=\max_i d_\eps(U_{x_i})<b$, a contradiction.
\end{proof}

\begin{lemma}\label{lem:source}
If $E\subset X$ is nonempty and compact and $\Dp(T,E)>a>0$, then
there is a countable compact set
\begin{equation}\label{eq:source-form}
 C=\{x_*\}\cup\{x_j:j\geq1\}\subset E,\qquad x_j\longrightarrow x_*,
\end{equation}
with distinct $x_j\ne x_*$ and $\Dp(T,C)>a$.
\end{lemma}

\begin{proof}
Choose $\eps_0>0$ and numbers $a<b_0<b_1<d_{\eps_0}(E)$.
Apply Lemma~\ref{lem:local-point} with $b=b_1$. For each $j\in\N$,
the set $U_j=E\cap B(x_*,1/j)$ has $d_{\eps_0}(U_j)\geq b_1>b_0$.
The definition of limsup therefore permits a strictly increasing
sequence $n_j$ and finite $(n_j,\eps_0)$-separated sets $F_j\subset U_j$
such that
\[
 \frac{\ell_p(\card{F_j})}{L_r(n_j)}>b_0.
\]
Put $C=\{x_*\}\cup\bigcup_jF_j$. Every neighborhood of $x_*$ contains
all but finitely many of these finite sets. Thus $C$ is compact and has
no accumulation point other than $x_*$. Since
$s_{n_j}(C,\eps_0)\geq\card{F_j}$, we have
$\Dp(T,C)\geq d_{\eps_0}(C)\geq b_0>a$. In particular $C$ is
infinite. Enumerating its distinct points other than $x_*$ gives
\eqref{eq:source-form}; every such enumeration converges to $x_*$,
because only finitely many points lie outside any neighborhood of $x_*$.
\end{proof}

\begin{lemma}\label{lem:tails}
Let $C$ have the form \eqref{eq:source-form}, and define
$C_J=\{x_*\}\cup\{x_j:j\geq J\}$. For every $L\in\N$ and $\rho>0$,
one has $C_J\subset B_L(x_*,\rho)$ for all sufficiently large $J$.
For each $J$ and $\eps>0$,
\[
 d_\eps(C_J)=d_\eps(C),\qquad \Dp(T,C_J)=\Dp(T,C).
\]
\end{lemma}

\begin{proof}
Each of the finitely many maps $T^i$, $0\leq i<L$, is continuous, so
$d_L(x_j,x_*)\to0$. This proves the first assertion. The other
assertions follow from finite deletion in Lemma~\ref{lem:finite}.
\end{proof}

\subsection{Trees from successive metric covers}

\begin{lemma}\label{lem:cover-tree}
Let $C\subset X$ be nonempty, let $L<N$ be positive integers, and let
$P\subset C\cap B_L(x_*,\rho)$ be a nonempty finite
$(N,\eta)$-separated set, where $\eta>2\rho$. Put $t_k=2^kL$ for
$1\leq k\leq J$, where $J$ is the first index with $t_J\geq N$.
There is a rooted tree on successive partitions of $P$ with these
properties:
\begin{enumerate}
\item Each level-$k$ atom has $d_{t_k}$-diameter less than $2\rho$.
\item If $B_k$ denotes the number of level-$k$ vertices, then
\begin{equation}\label{eq:node-product}
 B_k\leq\prod_{j=1}^k r_{t_j}(C,\rho)
 \leq r_{t_k}(C,\rho)^k.
\end{equation}
\item The terminal leaves are the individual points of $P$. If a
nonempty subset $F\subset P$ occupies $b_k$ vertices at level $k$, then
\begin{equation}\label{eq:occupied-cover}
 r_n(F,2\rho)\leq b_k\quad(n\leq t_k),\qquad
 r_n(F,2\rho)=1\quad(n\leq L).
\end{equation}
\end{enumerate}
\end{lemma}

\begin{proof}
At time $t_k$, take a cover of $C$ by exactly $r_{t_k}(C,\rho)$ open
Bowen balls of radius $\rho$ and order its members. Assign each point of
$P$ to the first ball containing it, discarding empty assigned sets.
This produces a partition $\mathcal P_k$ of $P$. Each atom is contained
in one ball and therefore has $d_{t_k}$-diameter less than $2\rho$.
Define the refining partitions
$\mathcal Q_k=\bigvee_{j=1}^k\mathcal P_j$, with $\mathcal Q_0=\{P\}$.
These define a rooted tree as in Section~\ref{sec:prelim}. The bound on
the number of intersections gives the first inequality in
\eqref{eq:node-product}; monotonicity in time gives the second.

An atom of $\mathcal Q_J$ cannot contain two distinct points of $P$:
their $d_N$-distance would be less than $2\rho<\eta$, contradicting
separation. Thus terminal atoms are singletons, and all vertices have
terminal descendants. For a selected subset $F$, choose one of its
points from each occupied level-$k$ atom. The radius-$2\rho$ balls
about those points cover $F$ at time $t_k$ and at every earlier time,
since $d_n\leq d_{t_k}$. Finally $F\subset B_L(x_*,\rho)$, so a single
ball centered at $x_*$ covers $F$ whenever $n\leq L$.
\end{proof}

\section{Exact extraction and hereditary lowerability}\label{sec:extraction}

The next construction adapts
\cite{Wang}*{Lemma~3.5}. Here the trees are supplied by
Lemma~\ref{lem:cover-tree}, and Lemma~\ref{lem:scale-calculus}
controls the resulting losses in the diagonal normalization.

\begin{proposition}\label{prop:extraction}
Let $p\geq2$ and let $C$ have the form \eqref{eq:source-form}. If
$H=\Dp(T,C)>a>0$, there is a countable compact set $K\subset C$ with
exactly one accumulation point such that $\Dp(T,K)=a$.
\end{proposition}

\begin{proof}
We use $r=p-1$ and the functions in \eqref{eq:fG}. By
Lemma~\ref{lem:range}, $0<H\leq1$.

\smallskip\noindent\emph{Choice of resolutions and parameters.}
Choose $\eps_m\downarrow0$ such that
\[
 h_m=d_{\eps_m}(C)\uparrow H,
\]
discarding finitely many initial terms to have $h_1>0$. Let $m(s)$ be
a sequence in which every positive integer appears infinitely often;
for instance, concatenate the finite lists
$(1),(1,2),(1,2,3),\ldots$. This recurrence will provide infinitely
many witness times at each fixed resolution before the resolution
limit is taken. Choose $\eta_s\downarrow0$ with
$0<\eta_s<h_1/2$, and set
\begin{equation}\label{eq:stage-parameters}
 c_s=h_{m(s)}-\eta_s,\qquad D_s=H+\eta_s.
\end{equation}
Then $0<c_s<H<D_s$ and $0<a<D_s$. Choose $\rho_s\downarrow0$ so that
$2\rho_s<\eps_{m(s)}$ for every $s$. The auxiliary number $D_s$ may
exceed one; it is a comparison parameter rather than a claimed entropy
value.

\smallskip\noindent\emph{A threshold before choosing a tail.}
Put $N_0=0$. At stage $s$, choose an integer $L_s>N_{s-1}$, in the
positive-logarithm domain, so large that
\begin{equation}\label{eq:stage-cover-upper}
 (1+\log_2 t)\log r_t(C,\rho_s)\leq f_{D_s}(t)
 \quad\text{for every integer }t\geq L_s.
\end{equation}
Here and below $\log_2 t=(\log t)/(\log2)$. To justify this choice,
take $H<D'_s<D_s$. Lemma~\ref{lem:critical} implies that the
fixed-resolution covering rate at $\rho_s$ is at most $H$, so
$r_t(C,\rho_s)\leq G_{D'_s}(t)$ for all sufficiently large integers
$t$. One can obtain this last inequality directly by inserting a
parameter strictly between $H$ and $D'_s$ in the proof of that lemma.
Equation \eqref{eq:absorb} then yields
$(1+\log_2t)f_{D'_s}(t)\leq f_{D_s}(t)$ eventually, establishing
\eqref{eq:stage-cover-upper}.

Only after $L_s$ has been fixed, choose a tail index $J_s$, strictly
larger than previous tail indices, for which
\begin{equation}\label{eq:stage-concentration}
 C_s:=\{x_*\}\cup\{x_j:j\geq J_s\}
 \subset B_{L_s}(x_*,\rho_s).
\end{equation}
This is possible by Lemma~\ref{lem:tails}. In particular, all times up
to $L_s$ are controlled by one ball.

\smallskip\noindent\emph{Separated source points and the selected cardinality.}
Finite deletion gives $d_{\eps_{m(s)}}(C_s)=h_{m(s)}>c_s$. There are
therefore arbitrarily large integers $N>L_s$ for which $C_s$ contains
an $(N,\eps_{m(s)})$-separated set $P$ with
\begin{equation}\label{eq:many-source-points}
 \card P\geq G_{c_s}(N).
\end{equation}
Indeed, choose $b$ with $c_s<b<h_{m(s)}$. Along arbitrarily large
times the normalized regularized logarithm of the separated count
exceeds $b$, while that of $G_{c_s}(N)$ tends to $c_s$.

For such times define
\[
 u_N=f_{D_s}^{-1}\left(\tfrac12 f_{c_s}(N)\right),\qquad
 M_N=\lfloor G_a(u_N)\rfloor.
\]
Lemma~\ref{lem:scale-calculus} allows us to choose one of the times
in \eqref{eq:many-source-points}, denoted by $N_s$, so large that
$u_{N_s}\geq L_s$, $u_{N_s}<N_s$, $M_{N_s}\geq2$, and
\begin{equation}\label{eq:stage-error}
 \left|\frac{\ell_p(M_{N_s})}{L_r(N_s)}-\frac{ac_s}{D_s}\right|
 <\frac1s.
\end{equation}
Write $u_s=u_{N_s}$, $M_s=M_{N_s}$, and let $P_s$ be the corresponding
separated set, with $A_s=\card{P_s}$. All these requirements are
eventual inequalities for the already fixed stage parameters, so they
can be met simultaneously along the unbounded available times.

Since $a<D_s$ and $L_r(u_s)>0$,
\begin{equation}\label{eq:Ms-upper}
 \log M_s\leq f_a(u_s)\leq f_{D_s}(u_s)
 =\tfrac12 f_{c_s}(N_s).
\end{equation}
Together with $A_s\geq\exp(f_{c_s}(N_s))$, this gives $1\leq M_s\leq A_s$.

\smallskip\noindent\emph{Selection on a metric covering tree.}
Apply Lemma~\ref{lem:cover-tree} to $P_s$, using $L=L_s$,
$N=N_s$, $\rho=\rho_s$, and $\eta=\eps_{m(s)}$. Its times are
$t_k=2^kL_s$, ending at the first $t_J\geq N_s$. Denote its node
counts by $B_{s,k}$. Since $k\leq1+\log_2t_k$, equations
\eqref{eq:node-product} and \eqref{eq:stage-cover-upper} give
\begin{equation}\label{eq:Bs-upper}
 B_{s,k}\leq r_{t_k}(C,\rho_s)^k\leq \exp(f_{D_s}(t_k)).
\end{equation}
Use Lemma~\ref{lem:tree-selection} to select $M_s$ leaves, and let
$F_s\subset P_s$ be the corresponding points. Write $b_{s,k}$ for the
occupied node counts. When $t_k\leq u_s$, we have
\begin{align*}
 \frac{M_sB_{s,k}}{A_s}
 &\leq\exp\bigl(f_a(u_s)+f_{D_s}(t_k)-f_{c_s}(N_s)\bigr)\\
 &\leq\exp\bigl(2f_{D_s}(u_s)-f_{c_s}(N_s)\bigr)=1.
\end{align*}
It follows that $b_{s,k}\leq4k^2$. When $t_k>u_s$, cardinality alone
gives $b_{s,k}\leq M_s\leq G_a(t_k)$. These two estimates hold for
every level in the same selected tree.

\smallskip\noindent\emph{One upper bound for every stage and every time.}
For $L_s<n\leq N_s$, let $t_k$ be the first dyadic time at least $n$.
Then $n\leq t_k<2n$ and $k\leq1+\log_2 n$. By
\eqref{eq:occupied-cover} and the two preceding estimates,
\[
 r_n(F_s,2\rho_s)\leq
 \max\{4(1+\log_2 n)^2,G_a(2n)\}.
\]
This includes the case $n<u_s<t_k$, for which the second estimate
applies. For $n\leq L_s$ the cover has one ball by
\eqref{eq:stage-concentration}. For $n>N_s$,
\[
 r_n(F_s,2\rho_s)\leq M_s\leq G_a(u_s)\leq G_a(n).
\]
Extend $G_a$ to a positive nondecreasing function on the remaining
bounded interval, keeping its original values for large arguments.
The three cases give a stage-independent bound
\begin{equation}\label{eq:Q-bound}
 \begin{gathered}
 r_n(F_s,2\rho_s)\leq Q(n)\qquad(n\geq1),\\
 Q(n)=4(1+\log_2(2n))^2(1+G_a(2n)).
 \end{gathered}
\end{equation}

\smallskip\noindent\emph{The compact union and its upper rate.}
Set
\begin{equation}\label{eq:K-union}
 K=\{x_*\}\cup\bigcup_{s\geq1}F_s.
\end{equation}
Each $F_s$ is finite and lies in $B(x_*,\rho_s)$ by
\eqref{eq:stage-concentration}. Since $\rho_s\to0$, the set $K$ is
countable and compact, and it has at most the accumulation point $x_*$.
Also $K\subset C$.

Fix $\eps>0$, and choose $s_0$ such that $2\rho_s<\eps$ whenever
$s\geq s_0$. Let $C_\eps=\sum_{s<s_0}\card{F_s}$. The integers
$L_s$ are strictly increasing, since
$L_s>N_{s-1}>L_{s-1}$. Hence at most $n$ stages satisfy $L_s<n$.
For the stages $s\geq s_0$ with $L_s<n$, use \eqref{eq:Q-bound}
and monotonicity in the radius. For every remaining stage $s\geq s_0$
with $L_s\geq n$, equation \eqref{eq:stage-concentration} implies
$F_s\subset B_n(x_*,\rho_s)\subset B_n(x_*,\eps)$. All those
remaining stages, together with $x_*$, need just one ball. Thus
\begin{equation}\label{eq:K-upper}
 r_n(K,\eps)\leq C_\eps+nQ(n)+1.
\end{equation}
For sufficiently large $n$, the right-hand side is at most
$C'_{\eps}n^3G_a(2n)$ for a fixed constant $C'_{\eps}>0$:
use $(1+\log_2(2n))^2=O(n^2)$ and $G_a(2n)\to\infty$.
Part (3) of Lemma~\ref{lem:scale-calculus} gives
\[
 \limsup_{n\to\infty}\frac{\ell_p(r_n(K,\eps))}{L_r(n)}\leq a.
\]
Since $\eps$ was arbitrary, $\Dp(T,K)\leq a$.

\smallskip\noindent\emph{A lower bound at each recurring resolution.}
Fix $m\in\N$. There are infinitely many stages with $m(s)=m$.
At each such stage, $F_s\subset K$ is $(N_s,\eps_m)$-separated,
so $s_{N_s}(K,\eps_m)\geq M_s$. Equation \eqref{eq:stage-error}
therefore implies
\[
 d_{\eps_m}(K)\geq
 \lim_{\substack{s\to\infty\\m(s)=m}}
 \frac{a(h_m-\eta_s)}{H+\eta_s}=\frac{ah_m}{H}.
\]
Letting $m\to\infty$ yields $\Dp(T,K)\geq a$. This proves equality.
Since $a>0$, Lemma~\ref{lem:finite} rules out a finite $K$; compactness
then ensures that $x_*$ is its unique accumulation point.
\end{proof}

\begin{proof}[Proof of Theorem~\ref{thm:main}]
If $a=0$, choose a singleton in $E$. If $a=\Dp(T,E)$, choose
$K_a=E$. For $0<a<\Dp(T,E)$, Lemma~\ref{lem:source} gives a
countable compact $C\subset E$ of the form \eqref{eq:source-form}
with $\Dp(T,C)>a$. Apply Proposition~\ref{prop:extraction} to $C$.
The resulting set lies in $E$ and has exactly one accumulation point.
\end{proof}

Thus, for every nonempty compact $E\subset X$ and each fixed $p\geq2$,
\[
 \{\Dp(T,K):\varnothing\ne K\subset E,\ K\text{ compact}\}
 =[0,\Dp(T,E)];
\]
the theorem gives inclusion from right to left, and monotonicity gives
the reverse inclusion.

\begin{remark}\label{rem:p1}
The proof uses $p\geq2$ at specific points. For a formal classical
replacement $f_b(t)=bt$, the estimate
$(\log t)f_b(t)=o(f_d(t))$ fails: the quotient equals
$(b/d)\log t$. In addition, $L_r(2n)/L_r(n)\to1$ for $r\geq1$,
whereas the classical denominator $n$ retains the factor two.
Finally, the source value $H$ is always finite here by
Lemma~\ref{lem:range}. Consequently this argument gives no hereditary
lowerability assertion for $h^{1,1}=h$ on infinite-entropy systems.
\end{remark}

\section{Sharpness and comparison with classical entropy}\label{sec:supplements}

We give two complements to Theorem~\ref{thm:main}. The first shows why
the assertion about one accumulation point is restricted to strictly
intermediate values. The second compares $\Dp$ with classical entropy
on the product-diagonal construction of Huang, Ye, and Zhang.

\subsection{The upper endpoint}\label{sec:endpoint}

\begin{proposition}\label{prop:endpoint}
For every integer $p\geq2$ and every $H\in(0,1]$, there are a compact
metric space $Z$, a homeomorphism $R\colon Z\to Z$, and a nonempty
compact set $E\subset Z$ such that
\[
 \Dp(R,E)=H,\qquad h(R,E)=0,
\]
while every nonempty compact subset $K\subset E$ with at most one
accumulation point satisfies $\Dp(R,K)<H$.
\end{proposition}

\begin{proof}
Choose $0<a_j<H$ with $a_j\uparrow H$. Since
$\Dp(\sigma,\Sigma_2)=1$, Theorem~\ref{thm:main} supplies nonempty
compact sets $A_j\subset\Sigma_2$ with $\Dp(\sigma,A_j)=a_j$.
Moreover, $h(\sigma,A_j)=0$: otherwise Lemma~\ref{lem:range} would
give $\Dp(\sigma,A_j)=1$, contrary to $a_j<1$.

Let $X_j$ be a separate copy of $\Sigma_2$, and set
$Z=\{\infty\}\cup\bigsqcup_{j\geq1}X_j$. Define
\[
 \begin{aligned}
 \delta(x,y)&=2^{-j}d_0(x,y) &&(x,y\in X_j),\\
 \delta(x,y)&=2^{-i}+2^{-j} &&(x\in X_i,\ y\in X_j,\ i\ne j),\\
 \delta(x,\infty)&=\delta(\infty,x)=2^{-j} &&(x\in X_j),
 \end{aligned}
\]
and $\delta(\infty,\infty)=0$. This is a metric. Within a component
the triangle inequality follows from that for $d_0$; a path leaving
$X_j$ has length at least $2\cdot2^{-j}$, whereas distances within
$X_j$ are at most $2^{-j}$. For points in different components the
displayed sums give the remaining triangle inequalities. Each $X_j$
is compact and open and closed in $Z$. A sequence with unbounded
component indices has a subsequence converging to $\infty$; otherwise
it has a subsequence in a fixed compact component or a constant
subsequence at $\infty$. Thus $Z$ is compact.

Let $R$ act as the shift on each $X_j$ and fix $\infty$. Both $R$ and
$R^{-1}$ are continuous within each component and preserve the distance
to $\infty$, so $R$ is a homeomorphism. Regard $A_j$ as a subset of
$X_j$ and put
\[
 E=\{\infty\}\cup\bigcup_{j\geq1}A_j.
\]
This set is closed: each $A_j$ is closed in its component, and limits
with unbounded component indices equal $\infty$. Hence $E$ is compact.

Use superscripts to distinguish covering numbers for $R$ and $\sigma$.
Fix $\eps>0$ and choose $J$ so that $2^{-j}<\eps$ for $j>J$.
All components with $j>J$ remain in the radius-$\eps$ Bowen ball
about $\infty$ at every time. Covering the other components separately
therefore gives
\begin{equation}\label{eq:endpoint-cover}
 r_n^R(E,\eps)\leq1+\sum_{j=1}^{J}
 r_n^\sigma(A_j,2^j\eps).
\end{equation}
A finite sum is at most its largest summand times the number of terms.
Thus, by the bounded-error estimate in Lemma~\ref{lem:finite}, its
normalized $\ell_p$-growth is at most the largest rate of its summands;
the same argument applies to the classical logarithm. It follows that
$h(R,E)=0$ and $\Dp(R,E)\leq H$ after taking the resolution limit.
For each fixed $j$ and $\eps<2^{-j}$, a ball meeting $X_j$ must have
its center in $X_j$, already by the distance at time zero. Hence
\[
 r_n^R(A_j,\eps)=r_n^\sigma(A_j,2^j\eps)
 \qquad(\eps<2^{-j}).
\]
Consequently $\Dp(R,A_j)=a_j$, and monotonicity yields
$\Dp(R,E)\geq\sup_j a_j=H$.

Now let $K\subset E$ be nonempty and compact with at most one
accumulation point. If there is none, $K$ is finite and $\Dp(R,K)=0$.
If the only accumulation point is $\infty$, each $K\cap A_j$ is
finite: it is compact, and an infinite such intersection would have
an accumulation point in $X_j$. At each fixed resolution only finitely
many components require separate covers in \eqref{eq:endpoint-cover};
their intersections with $K$ are finite, and the remaining points fit
in one ball about $\infty$. Thus $r_n^R(K,\eps)$ is bounded in $n$
and $\Dp(R,K)=0$.

Finally, suppose the unique accumulation point lies in $A_j$. Then
$K$ meets only finitely many components, since otherwise a sequence
with increasing component indices would have the second accumulation
point $\infty$. Its intersections with the other components are
finite by compactness. Therefore $K\setminus A_j$ is finite, and
Lemma~\ref{lem:finite} gives $\Dp(R,K)\leq a_j<H$.
\end{proof}

\subsection{The product-diagonal comparison}\label{sec:classical}

Let $(X,T)$ be a compact homeomorphism system with
$0<h(T,X)<\infty$, give $X$ a compatible metric $d\leq1$, and set
\[
 Y=X^{\N},\qquad S(y)_j=T^j y_j,\qquad
 \Delta x=(x,x,\ldots).
\]
Equip $Y$ with the metric \eqref{eq:product-metric}. The maps $S$ and
$S^{-1}$ are continuous coordinatewise, so $S$ is a homeomorphism.
The map $\Delta$ is a homeomorphism onto its image, with inverse given
by the first coordinate projection. In this construction,
\cite{HYZ14}*{Proposition~3.2 and its proof} gives, for every nonempty
compact $A\subset X$,
\begin{equation}\label{eq:diagonal-classical}
 h(S,\Delta A)=
 \begin{cases}
 0,&h(T,A)=0,\\
 +\infty,&h(T,A)>0.
 \end{cases}
\end{equation}
In particular, $\Delta X$ has infinite classical entropy and no compact
subset of finite positive classical entropy. The diagonal scales have
the following behavior.

\begin{proposition}\label{prop:diagonal-entropy}
For every nonempty compact $A\subset X$ and every $p\geq2$,
\begin{equation}\label{eq:diagonal-scale}
 \Dp(S,\Delta A)=\Dp(T,A).
\end{equation}
\end{proposition}

\begin{proof}
Projection onto the first coordinate sends a radius-$\eps$ Bowen ball
for $S$ into a radius-$2\eps$ Bowen ball for $T$, since
$\rho(S^iy,S^iz)\geq\frac12d(T^iy_1,T^iz_1)$ for every $i$.
The centers project to points of $X$, whether or not they belong to $A$.
Projecting a cover therefore yields
\begin{equation}\label{eq:diagonal-lower}
 r_n^T(A,2\eps)\leq r_n^S(\Delta A,\eps),
\end{equation}
and hence $\Dp(T,A)\leq\Dp(S,\Delta A)$.

For the reverse inequality, fix $\eps>0$, choose $m$ with
$\sum_{j>m}2^{-j}<\eps/2$, and take $0<\delta<\eps/2$.
If $d_{mn}^T(x,u)<\delta$, then for $0\leq i<n$,
\[
 \rho(S^i\Delta x,S^i\Delta u)
 \leq\sum_{j=1}^{m}2^{-j}d(T^{ji}x,T^{ji}u)
      +\sum_{j>m}2^{-j}
 <\delta+\eps/2<\eps,
\]
because $ji<mn$ for $j\leq m$. Replacing every center $u$ in a
length-$mn$ cover of $A$ by $\Delta u$ gives
\begin{equation}\label{eq:diagonal-upper}
 r_n^S(\Delta A,\eps)\leq r_{mn}^T(A,\delta).
\end{equation}
Put $r=p-1\geq1$. The integer $m$ is fixed at this resolution, and
$L_r(mn)/L_r(n)\to1$. Thus
\begin{align*}
 \limsup_{n\to\infty}
 \frac{\ell_p(r_n^S(\Delta A,\eps))}{L_r(n)}
 &\leq\limsup_{n\to\infty}
 \frac{\ell_p(r_{mn}^T(A,\delta))}{L_r(mn)}\\
 &\leq\limsup_{N\to\infty}
 \frac{\ell_p(r_N^T(A,\delta))}{L_r(N)}
 \leq\Dp(T,A).
\end{align*}
The middle inequality only bounds a subsequence upper limit by the
full upper limit. Letting $\eps\downarrow0$ proves the reverse bound.
\end{proof}

For the binary full shift, Theorem~\ref{thm:main} supplies, for each
fixed $p\geq2$ and $0<a<1$, a nonempty compact set $A\subset\Sigma_2$
with $\Dp(\sigma,A)=a$. Lemma~\ref{lem:range} then gives
$h(\sigma,A)=0$. Equations \eqref{eq:diagonal-classical} and
\eqref{eq:diagonal-scale} imply
\[
 h(S,\Delta A)=0,\qquad \Dp(S,\Delta A)=a.
\]
Thus the classical entropy restriction on $\Delta\Sigma_2$ allows
intermediate diagonal scale values on its classical zero-entropy subsets.



\begin{bibdiv}
\begin{biblist}

\bib{AKM}{article}{
      author={Adler, R.~L.},
      author={Konheim, A.~G.},
      author={McAndrew, M.~H.},
       title={Topological entropy},
        date={1965},
     journal={Trans. Amer. Math. Soc.},
      volume={114},
      number={2},
       pages={309\ndash 319},
  note={\href{https://doi.org/10.1090/S0002-9947-1965-0175106-9}{doi:\nolinkurl{10.1090/S0002-9947-1965-0175106-9}}},
}

\bib{Bowen71}{article}{
      author={Bowen, R.},
       title={Entropy for group endomorphisms and homogeneous spaces},
        date={1971},
     journal={Trans. Amer. Math. Soc.},
      volume={153},
       pages={401\ndash 414},
  note={\href{https://doi.org/10.1090/S0002-9947-1971-0274707-X}{doi:\nolinkurl{10.1090/S0002-9947-1971-0274707-X}}},
}

\bib{Bowen73}{article}{
      author={Bowen, R.},
       title={Topological entropy for noncompact sets},
        date={1973},
     journal={Trans. Amer. Math. Soc.},
      volume={184},
       pages={125\ndash 136},
  note={\href{https://doi.org/10.1090/S0002-9947-1973-0338317-X}{doi:\nolinkurl{10.1090/S0002-9947-1973-0338317-X}}},
}

\bib{Carvalho}{article}{
      author={de~Carvalho, M.},
       title={Entropy dimension of dynamical systems},
        date={1997},
     journal={Portugal. Math.},
      volume={54},
      number={1},
       pages={19\ndash 40},
         url={https://eudml.org/doc/48825},
}

\bib{DHP}{article}{
      author={Dou, D.},
      author={Huang, W.},
      author={Park, K.~K.},
       title={Entropy dimension of topological dynamical systems},
        date={2011},
     journal={Trans. Amer. Math. Soc.},
      volume={363},
      number={2},
       pages={659\ndash 680},
  note={\href{https://doi.org/10.1090/S0002-9947-2010-04906-2}{doi:\nolinkurl{10.1090/S0002-9947-2010-04906-2}}},
}

\bib{Helfter}{article}{
      author={Helfter, M.},
       title={Scales},
        date={2025},
     journal={Math. Z.},
      volume={310},
      number={1},
       pages={Paper No.~15},
  note={\href{https://doi.org/10.1007/s00209-025-03719-5}{doi:\nolinkurl{10.1007/s00209-025-03719-5}}},
}

\bib{HYZ10}{article}{
      author={Huang, W.},
      author={Ye, X.},
      author={Zhang, G.},
       title={Lowering topological entropy over subsets},
        date={2010},
     journal={Ergodic Theory Dynam. Systems},
      volume={30},
      number={1},
       pages={181\ndash 209},
  note={\href{https://doi.org/10.1017/S0143385709000066}{doi:\nolinkurl{10.1017/S0143385709000066}}},
}

\bib{HYZ14}{article}{
      author={Huang, W.},
      author={Ye, X.},
      author={Zhang, G.},
       title={Lowering topological entropy over subsets revisited},
        date={2014},
     journal={Trans. Amer. Math. Soc.},
      volume={366},
      number={8},
       pages={4423\ndash 4442},
  note={\href{https://doi.org/10.1090/S0002-9947-2014-06117-5}{doi:\nolinkurl{10.1090/S0002-9947-2014-06117-5}}},
}

\bib{KT}{article}{
      author={Katok, A.},
      author={Thouvenot, J.-P.},
       title={Slow entropy type invariants and smooth realization of commuting
  measure-preserving transformations},
        date={1997},
     journal={Ann. Inst. H. Poincar{\'e} Probab. Statist.},
      volume={33},
      number={3},
       pages={323\ndash 338},
  note={\href{https://doi.org/10.1016/S0246-0203(97)80094-5}{doi:\nolinkurl{10.1016/S0246-0203(97)80094-5}}},
}

\bib{Wang}{misc}{
      author={Wang, X.},
       title={Hereditary lowerability of topological dynamical systems},
        date={2026},
         url={https://arxiv.org/abs/2608.07115v1},
        note={Preprint,
  \href{https://arxiv.org/abs/2608.07115v1}{arXiv:2608.07115v1}},
}

\bib{XCZ}{misc}{
      author={Xie, Y.},
      author={Chen, E.},
      author={Zhou, X.},
       title={Entropy scales in topological dynamical systems},
        date={2026},
         url={https://arxiv.org/abs/2608.10460v2},
        note={Preprint,
  \href{https://arxiv.org/abs/2608.10460v2}{arXiv:2608.10460v2}},
}

\end{biblist}
\end{bibdiv}
\end{document}